\documentclass{svmult}

\usepackage{graphicx}        
\usepackage{multicol}        
\usepackage[bottom]{footmisc}

\usepackage[latin9]{inputenc}
\usepackage{mathtools}
\usepackage[english]{babel}
\usepackage{amsfonts}

\usepackage{mathptmx}       
\usepackage{helvet}         
\usepackage{courier}        
\usepackage{type1cm}        

\usepackage[a4paper]{geometry}

\renewcommand{\hom}{\mathrm{Hom}}
\def\quot{/\!\!/}

\newcommand{\pexp}{\operatorname{PExp}}

\begin{document}

\title*{On Hodge polynomials of Singular Character Varieties}
\author{Carlos Florentino\inst{1}\and
Azizeh Nozad\inst{2}\and Jaime Silva\inst{3}\and Alfonso Zamora\inst{4}}

\institute{Departamento de Matem\'{a}tica, Faculdade de Ci\^encias, Univ. de Lisboa, Campo Grande, Edf. C6, Lisbon, Portugal 
\texttt{caflorentino@ciencias.ulisboa.pt}
\and School of Mathematics, Institute for Research in Fundamental Sciences (IPM), P.O.Box: 19395-5746, Tehran, Iran 
\texttt{anozad@ipm.ir}
\and Centro de Matem\'atica da Universidade do Porto, Faculdade de Ci\^encias da Universidade do Porto, Rua do Campo Alegre, 4169-007 Porto, Portugal\texttt{jaime.silva@fc.up.pt}
\and Departamento de Matem\'atica Aplicada a las TIC, ETSI Inform\'aticos, Universidad Polit\'ecnica de Madrid, Campus de Montegancedo,
28660, Madrid, Spain
\texttt{alfonso.zamora@upm.es}}

%
%
\maketitle

\subsection*{Abstract}
Let $\mathcal{X}_{\Gamma}G:=\hom(\Gamma,G)\quot G$ be the $G$-character
variety of $\Gamma$, where $G$ is a complex reductive group and
$\Gamma$ a finitely presented group. We introduce new techniques
for computing Hodge-Deligne and Serre polynomials of $\mathcal{X}_{\Gamma}G$,
and present some applications, focusing on the cases when $\Gamma$
is a free or free abelian group. Detailed constructions and proofs
of the main results will appear elsewhere.

\section{Introduction}

Let $G$ be a connected reductive complex algebraic group, and $\Gamma$
be a finitely presented group. The $G$-character variety of $\Gamma$
is defined to be the (affine) geometric invariant theory (GIT) quotient
\[
\mathcal{X}_{\Gamma}G=\hom(\Gamma,G)\quot G.
\]
The most well studied families of character varieties include the
cases when the group $\Gamma$ is the fundamental group of a Riemann
surface $\Sigma$, and its ``twisted'' variants. In these cases,
the non-abelian Hodge correspondence (see, for example \cite{Si})
shows that (components of) $\mathcal{X}_{\Gamma}G$ are homeomorphic
to certain moduli spaces of $G$-Higgs bundles which appear in connection
to important problems in Mathematical-Physics: for example, these
spaces play an important role in the quantum field theory interpretation
of the geometric Langlands correspondence, in the context of mirror
symmetry (\cite{KW}).

The study of geometric and topological properties of character varieties
is an active topic and there are many recent advances in the computation
of their Poincar\'e polynomials and other invariants. For the surface
group case ($\Gamma=\pi_{1}(\Sigma)$ and related groups) the calculations
of Poincar\'e polynomials started with Hitchin and Gothen, and have
been pursued more recently by Hausel, Lettelier, Mellit, Rodriguez-Villegas,
Schiffmann and others, who also considered the parabolic version of
these character varieties (see \cite{HRV,Me,Sc}). Those recent results
use arithmetic methods: it is shown that the number of points of the
corresponding moduli space over finite fields is given by a polynomial,
which turns out to coincide with the $E$-polynomial of $\mathcal{X}_{\Gamma}G$
(\cite[Appendix]{HRV}). Then, in the smooth case, the pure nature
of the cohomology of Higgs bundles moduli spaces allows the derivation
of the Poincar\'e polyomial from the $E$-polynomial.

On the other hand, for many important classes of \emph{singular} character
varieties, explicitly computable formulas for the $E$-polynomials
(also called Serre, or Hodge-Euler polynomials) are very hard to obtain. In
the articles of Logares, Mu\~noz, Newstead and Lawton \cite{LMN}, \cite{LM}
(using geometric methods) and of Baraglia and Hekmati \cite{BH} (using
arithmetic methods), the $E$-polynomials are computed for several
character varieties, with $G=GL(n,\mathbb{C})$, $SL(n,\mathbb{C})$
and $PGL(n,\mathbb{C})$ for small values of $n$, but the computations
quickly become intractable for $n$ higher than $3$.

In this short article, we describe some of the techniques and constructions
that we have recently developed for computations of $E$-polynomials
of singular character varieties, and present some of their main applications. 

The outline of the article is as follows. Section 2 covers notations
and preliminaries on mixed Hodge and $E$-polynomials and on character
varieties in the context of GIT. In section 3, we explain how to use
equivariant mixed Hodge structures to study (the identity component
of) $\mathcal{X}_{\Gamma}G$ when $\Gamma$ is a free abelian group
and $G$ a classical group. These character vareities have orbifold
singularities and we can obtain their full mixed Hodge polynomials.
In section 4, for arbitrary $\Gamma$, we define a stratification
of $GL(n,\mathbb{C})$-character varieties (which also exists for
$G=SL(n,\mathbb{C})$ or $PGL(n,\mathbb{C})$) which allows writing
down an explicit plethystic exponential relation between generating
functions of the $E$-polynomials of $\mathcal{X}_{\Gamma}GL(n,\mathbb{C})$
and of its locus of irreducible representations $\mathcal{X}_{\Gamma}^{irr}GL(n,\mathbb{C})$.
Finally, in Section 5, we consider the free group $\Gamma=F_{r}$
of rank $r$, and announce the solution of a conjecture of Lawton
and Mu\~noz: the $E$-polynomials of $\mathcal{X}_{F_{r}}SL(n,\mathbb{C})$
and of $\mathcal{X}_{F_{r}}PGL(n,\mathbb{C})$ coincide, for every
$n\in\mathbb{N}$. For lack of space, the proofs are omitted and will
be published elsewhere.

\subsection*{Acknowledgements }

We would like to thank S. Lawton, M. Logares, S. Mozgovoy, V. Mu\~noz,
A. Oliveira and F. Rodriguez-Villegas for several interesting and
very useful conversations on topics around mixed Hodge structures
and $E$-polynomials. We also thank Alexander Schmitt for the invitation
to deliver a talk at the Special Session on Complex Geometry of the
conference ISAAC 2019, where these results were presented.

This work was partially supported by CAMGSD and CMAFcIO of the University of Lisbon, the projects PTDC/MAT-GEO/2823/2014 and PTDC/MAT-PUR/30234/2017, FCT Portugal, a grant from IPM (Iran) and ICTP (Italy), and the project MTM2016-79400-P by the Spanish government.

\section{Preliminaries on Hodge-Deligne polynomials, Affine GIT and Character
Varieties}

In this article, all algebraic varieties are defined over $\mathbb{C}$,
$G$ is a connected reductive algebraic group, and $\Gamma$ is a finitely presented group.


Let $X$ be a quasi-projective variety (not necessarily irreducible), 
of complex dimension $\leq d$.
Deligne showed that the compactly supported
cohomology $H_{c}^{*}(X):=H_{c}^{*}(X,\mathbb{C})$ can be endowed
with a mixed Hodge structure whose mixed Hodge numbers are given by
\[
h^{k,p,q}(X):=\dim_{\mathbb{C}}H_{c}^{k,p,q}(X)\in\mathbb{N}_{0},
\]
for $k,p,q\in\{0,\cdots,2d\}$, and we call $(p,q)$ the $k$-weights
of $X$, if $h^{k,p,q}\neq0$ (c.f. \cite{De}, \cite{PS}). 

Mixed Hodge numbers are symmetric in the weights, $h^{k,p,q}=h^{k,q,p}$,
and $\dim_{\mathbb{C}}H_{c}^{k}(X)=\sum_{p,q}h^{k,p,q}$. Therefore,
they provide the (compactly supported) Betti numbers, yielding the
usual Betti numbers, by Poincar\'e duality, in the non-singular case. They
are also the coefficients of the mixed Hodge polynomial of $X$ on
three variables, 
\begin{equation}
\mu(X;\,t,u,v):=\sum_{k,p,q}h^{k,p,q}(X)\ t^{k}u^{p}v^{q}\in\mathbb{\mathbb{N}}_{0}[t,u,v],\label{eq:mu}
\end{equation}
which specializes to the (compactly supported) Poincar\'e polynomial
by setting $u=v=1$, $P_{t}^{c}(X):=\mu(X;\,t,1,1)$ (and provides
the usual Poincar\'e polynomial in the smooth situation). Plugging $t=-1$,
mixed Hodge polynomials convert into the $E$-polynomial of $X$,
or the \emph{Serre polynomial} of $X$, given by 
\[
E(X;\,u,v)=\sum_{k,p,q}(-1)^{k}h^{k,p,q}(X)\ u^{p}v^{q}\in\mathbb{Z}[u,v].
\]
From the $E$-polynomial we can compute the (compactly supported)
Euler characteristic of $X$ as $\chi^{c}(X)=E(X;\,1,1)=\mu(X;\,-1,1,1)$.

Serre polynomials satisfy an \emph{additive property with respect
to stratifications} by locally closed (in the Zariski topology) strata:
if $X$ has a closed subvariety $Z\subset X$ we have (see, eg. \cite{PS}),
\[
E(X)=E(Z)+E(X\setminus Z).
\]
The $E$-polynomial also satisfies (c.f. \cite{DL,LMN}) a \emph{multiplicative
property for fibrations. }Namely, for a given algebraic fibration
$F\hookrightarrow X\to B$, we have 
\[
E(X)=E(F)\cdot E(B)
\]
in any of the following three situations:
\begin{itemize}
\item[(i)]  the fibration is locally trivial in the Zariski topology of $B$, 
\item[(ii)]  $\,$ $F$, $X$ and $B$ are smooth, the fibration is locally trivial
in the complex analytic topology, and $\pi_{1}(B)$ acts trivially
on $H_{c}^{*}(F)$, or 
\item[(iii)] $\, \,$ $X$, $B$ are smooth and $F$ is a complex connected Lie group.
\end{itemize}
We say $X$ is of \emph{Hodge-Tate type} (also called \emph{balanced}
type) if all the $k$-weights are of the form $(p,p)$ with $p\in\{0,\cdots,k\}$,
in which case the sum in $\mu(X)$ reduces to a one-variable sum.
In particular, the $E$-polynomials of Hodge-Tate type varieties depend
only the product $uv$, so we write $x=uv$ and use the notation $E(X;\,x):=E(X;\,\sqrt{x},\sqrt{x})\in\mathbb{Z}[x].$


Now let $X$ be an affine algebraic variety, and
let the reductive group $G$ act algebraically on $X$. The induced action of $G$ on
the ring $\mathbb{C}[X]$ of regular functions on $X$ defines the
(affine) GIT quotient 
\[
X\quot G:=Spec\left(\mathbb{C}[X]^{G}\right),
\]
where $\mathbb{C}[X]^{G}$ is the subring of $G$-invariants in $\mathbb{C}[X]$.
This quotient identifies $G$-orbits whose closures intersect, such
that each point in the quotient classifies an equivalence class of
orbits, leading to a stability condition. Let $G_{x}\subset G$ be
the stabilizer of $x\in X$ and consider the orbit map through $x$,
$\psi_{x}:G\to X;\;g\mapsto g\cdot x$. We define $x\in X$ to be
\emph{stable} if $\psi_{x}$ is a proper map and \emph{polystable}
if the orbit $G\cdot x$ is closed in $X$. Stability implies polystability,
but not conversely.

GIT shows that the stable locus $X^{s}\subset X$ is a Zariski open
set (hence dense, when non-empty) and that the restriction of the
affine quotient map $\Phi:X\to X\quot G$ to the stable locus, $X^{s}\to X^{s}/G$,
is a geometric quotient (or an orbit space), where $\Phi(X^{s})$
is Zariski open in $X\quot G$.

Now, consider a finitely presented group $\Gamma$. The (generally
singular) algebraic variety of representations of $\Gamma$ in $G$
is 
\[
\mathcal{R}_{\Gamma}G=\hom(\Gamma,G).
\]
Each $\rho\in\mathcal{R}_{\Gamma}G$ is determined by $\rho(\gamma)$,
for each generator $\gamma\in\Gamma$, and satisfying the relations
of the group $\Gamma$. There is an algebraic action of $G$ on the
variety $\mathcal{R}_{\Gamma}G$ by conjugation of representations,
$g^{-1}\rho g$, yielding the $G$-character variety of $\Gamma$,
\[
\mathcal{X}_{\Gamma}G:=\hom(\Gamma,G)\quot G,
\]
as the GIT quotient.

By definition, polystable representations are representations $\rho:\Gamma\to G$
whose orbits $G\cdot\rho:=\{g\rho g^{-1}:\ g\in G\}$ are Zariski
closed in $\mathcal{R}_{\Gamma}G$. Alternatively, a representation
$\rho$ is polystable if and only if it is completely reducible (i.e,
if $\rho(\Gamma)\subset P\subset G$ for some proper parabolic $P$
of $G$, then $\rho(\Gamma)$ is contained in a Levi subgroup of $P$). 
Denote the subset of polystable representations in
$\mathcal{R}_{\Gamma}G$ by $\mathcal{R}_{\Gamma}^{ps}G\subset\mathcal{R}_{\Gamma}G$,
which is a Zariski locally-closed subvariety containing the stable
locus $\mathcal{R}_{\Gamma}^{s}G\subset\mathcal{R}_{\Gamma}G$. 
\begin{proposition}

\cite{FL1}\label{prop:character-polystable representation} There
is a bijective correspondence: 
\[
\mathcal{X}_{\Gamma}G=\mathcal{R}_{\Gamma}G\quot G\cong\mathcal{R}_{\Gamma}^{ps}G/G,
\]
where the right hand side is called the polystable quotient. 
\end{proposition}

We say that $\rho$ is \emph{irreducible} if $\rho(\Gamma)$
is not contained in a proper parabolic subgroup of $G$. 
Alternatively, $\rho$ is irreducible  if it is polystable and
$Z_{\rho}$, the centralizer of $\rho(\Gamma)$ inside $G$, is a finite 
extension of the center $ZG\subset G$. 
Denote by $\mathcal{R}_{\Gamma}^{irr}G\subset\mathcal{R}_{\Gamma}^{ps}G$
the subset of irreducible representations (being a Zariski open subset
of $\mathcal{R}_{\Gamma}G$, $\mathcal{R}_{\Gamma}^{irr}G$ is a quasi-projective
variety), and since irreducibility is well defined on $G$-orbits,
denote by 
\begin{equation}
\mathcal{X}_{\Gamma}^{irr}G:=\mathcal{R}_{\Gamma}^{irr}G/G\label{eq:X-irr}
\end{equation}
the\emph{ $G$-irreducible character variety of $\Gamma$}, which
is a geometric quotient, as it happens with the stable locus. In fact,
it can be proved that irreducibility is equivalent to GIT stability
for character varieties (see \cite[Thm. 1.3(1)]{CF}).

\section{The Free Abelian Case}

In this section, we are concerned with the determination of the mixed
Hodge polynomials of character varieties $\mathcal{X}_{\Gamma}G$
of the free abelian group of rank $r$, $\Gamma\cong\mathbb{Z}^{r}$. 
As we always work over $\mathbb{C}$, we abbreviate the notation of
the classical groups such as the linear group, special linear, special
orthogonal and symplectic to $GL_{n}$, $SL_{n}$, $SO_{n}$ and $Sp_{n}$,
respectively (instead of $GL(n,\mathbb{C})$, etc).


The topology and geometry of the character varieties $\mathcal{X}_{\mathbb{Z}^{r}}G$
was studied in \cite{FL2,Sk}, among others. Most important for
us are the following facts:
\begin{itemize}
\item[(i)]  there is only one irreducible component containing the trivial representation,
that we denote by $\mathcal{X}_{\mathbb{Z}^{r}}^{0}G$ \cite[Theorem 2.1]{Sk},
\item[(ii)] $\,$ if the semisimple part of $G$ is a classical group (ie, one of
$SL_{n}$, $SO_{n}$ and $Sp_{n}$), there exists an algebraic isomorphism
\begin{eqnarray}
\mathcal{X}_{\mathbb{Z}^{r}}^{0}G & \cong & \left(T_{G}\right)^{r}/W_{G}\label{eq:algis}
\end{eqnarray}
where $T_{G}$ is a maximal torus of $G$, and $W_{G}$ its Weyl group
\cite[Theorem 2.1]{Sk},
\item[(iii)] $\, \,$ the irreducibility of the free abelian character varieties $\mathcal{X}_{\mathbb{Z}^{r}}G$
can be characterized, in terms of $G$: for example, if the semisimple
part of $G$ is a product of $SL_{n}$'s and $Sp_{n}$'s then $\mathcal{X}_{\mathbb{Z}^{r}}G$
is irreducible, so that $\mathcal{X}_{\mathbb{Z}^{r}}G\cong\mathcal{X}_{\mathbb{Z}^{r}}^{0}G$
\cite[Theorem 1.2]{FL2}.
\end{itemize}
We now focus on the determination of the mixed Hodge numbers of $\mathcal{X}_{\mathbb{Z}^{r}}G$
when it is irreducible, or of $\mathcal{X}_{\mathbb{Z}^{r}}^{0}G$
when the algebraic isomorphism \eqref{eq:algis} applies. We start
by explaining how mixed Hodge numbers transform under finite quotients.


Let $X$ be a complex quasi-projective variety and $F$ a finite group
acting algebraically on it. The action of $F$ on $X$ induces an
action on its cohomology. Since $F$ acts by algebraic isomorphisms,
it also induces an action on the mixed Hodge components. Then we can
regard $H^{k,p,q}\left(X\right)$ as $F$-modules, that we denote
by $\left[H^{k,p,q}\left(X\right)\right]_{F}$. As in equation \eqref{eq:mu}
for the mixed Hodge polynomial, we codify these in the \emph{equivariant
mixed Hodge polynomial}, defined by 
\begin{eqnarray*}
\mu_{F}(X;\,t,u,v) & \coloneqq & \sum_{k,p,q}\left[H^{k,p,q}\left(X\right)\right] 
t^{k}u^{p}v^{q}\in R\left(F\right)\left[t,u,v\right]
\end{eqnarray*}
whose coefficients belong to $R(F)$, the representation ring of $F$.
The polynomial $\mu_{F}(X;\,t,u,v)$ may also be seen as a polynomial-weighted
representation. For instance, one can consider equivariant
cohomology to obtain an isomorphism 
\begin{eqnarray}
H^{*}\left(X/F\right) & \cong & H^{*}\left(X\right)^{F}\label{eq:quoinv}
\end{eqnarray}
that respects mixed Hodge structures. In particular,
this isomorphism allows us to identify the mixed Hodge polynomial
of the quotient $X/F$ as the coefficient of the trivial representation
of $\mu_{F}(X;\,t,u,v)$ when written on a basis of irreducible representations
of $F$. Another important consequence for us is the inequality 
$h^{k,p,q}\left(X\right)\geq h^{k,p,q}\left(X/F\right)$,
which holds since $H^{k,p,q}\left(X/F\right)$ is given by the $F$-invariant
part of $H^{k,p,q}\left(X\right)$. We conclude that if $X$ is, for
instance, a balanced variety, or if its mixed Hodge structure is actually
\emph{pure} (that is, if $h^{k,p,q}\neq0$ then $k=p+q$), then the same
holds for $X/F$. 

We now summarize our strategy to obtain the mixed Hodge polynomials
of $\mathcal{X}_{\mathbb{Z}^{r}}^{0}G$, in the cases when the isomorphism
\eqref{eq:algis} holds (so, these character varieties are
isomorphic to finite quotients of algebraic tori).
The only non-zero Hodge numbers of the maximal torus $T_{G}\cong\left(\mathbb{C}^{*}\right)^{n}$
are $h^{k,k,k}\left(T_{G}\right)$. Moreover, its natural mixed Hodge structure satisfies:
\begin{eqnarray*}
H^{k,k,k}\left(T_{G}\right) & \cong & \bigwedge^{k} \, \, H^{1,1,1}\left(T_{G}\right).
\end{eqnarray*}
So, the action of $W_G$ on the cohomology ring can be understood from the 
one on the mixed Hodge component $H^{1,1,1}\left(T_{G}\right)$. 
The next three theorems are proved
in \cite{FS}.
\begin{theorem}
\label{thm:det-formula} For a reductive group $G$ satisfying \eqref{eq:algis},
we have 
\begin{eqnarray*}
\mu(\mathcal{X}_{\mathbb{Z}^{r}}^{0}G;\,t,u,v) & = & \frac{1}{|W_{G}|}\sum_{g\in W_{G}}\left[\det\left(I+tuv\,A_{g}\right)\right]^{r}
\end{eqnarray*}
where $A_{g}$ is the automorphism of $H^{1,1,1}(T_{G})$ induced
by the action of $g\in W_{G}$. 
\end{theorem}
The proof starts by establishing the $r=1$ case, and using the diagonal
action for higher $r$ as well as the isomorphism \eqref{eq:algis},
together with the multiplicative relation for the equivariant polynomials
$\mu_{W_{G}}\left(T_{G}^{r}\right)=\mu_{W_{G}}\left(T_{G}\right)^{\otimes r}.$
We remark that Theorem \ref{thm:det-formula} generalizes a formula for the Poincar\'e
polynomial of $\mathcal{X}_{\mathbb{Z}^{r}}^{0}G$, recently obtained
in \cite{St}.


To further work with Theorem \ref{thm:det-formula}, we examine the
induced action of $W_{G}$ on $H^{1,1,1}\left(T_{G}\right)$ for some
classical groups. In the case $G=GL_{n}$, the Weyl group is the symmetric
group $S_{n}$ on $n$ letters, which acts on $H^{1,1,1}\left(T_{G}\right)\cong\mathbb{C}^{n}$
by permutation of coordinates, and we obtain a general formula in
terms of partitions of $n$. 

A partition of $n\in\mathbb{N}$ is denoted by $[k]=[1^{k_{1}}\cdots j^{k_{j}}\cdots n^{k_{n}}]$
where the exponent $k_{j}\geq0$ is the number of parts of size $j\in\{1,\cdots,n\}$,
so that $n=\sum_{j=1}^{n}j\cdot k_{j}$. Let $\mathcal{P}_{n}$ denote
the finite set of partitions of $n$.
\begin{theorem}
The mixed Hodge polynomials of $\mathcal{X}_{\mathbb{Z}^{r}}GL_{n}$
and of $\mathcal{X}_{\mathbb{Z}^{r}}SL_{n}$ satisfy 
\[
\mu(\mathcal{X}_{\mathbb{Z}^{r}}GL_{n};\,t,u,v)\ =\ \mu(\mathcal{X}_{\mathbb{Z}^{r}}SL_{n};\,t,x)(1+tuv)^{r}=\sum_{[k]\in\mathcal{P}_{n}}\prod_{j=1}^{n}\frac{(1-(-tuv)^{j})^{k_{j}r}}{k_{j}!\,j^{k_{j}}},
\]

\end{theorem}
By using similar considerations as for the $GL_{n}$ case, we can
also deduce a concrete formula for $Sp_{n}$ in terms of \emph{bipartitions}.
A bipartition of $n$, denoted $[a,b]\in\mathcal{B}_{n}$ consists
of two partitions $[a]\in\mathcal{P}_{k}$ and $[b]\in\mathcal{P}_{l}$,
such that $0\leq k,l\leq n$ with $k+l=n$.
One can show that bipartitions of $n$ are in one-to-one correspondence
with conjugacy classes in $W_{Sp_{n}}$, the Weyl group of $Sp_{n}$. 
\begin{theorem}
The mixed Hodge polynomial of $\mathcal{X}_{\mathbb{Z}^{r}}Sp_{n}\mathbb{C}$
is given by
\begin{eqnarray*}
\mu(\mathcal{X}_{\mathbb{Z}^{r}}Sp_{n};\,t,u,v) & = & \frac{1}{2^{n}n!}\sum_{[a,b]\in\mathcal{B}_{n}}c_{[a,b]}\prod_{i=1}^{k}(1-(-tuv)^{i})^{a_{i}r}\prod_{j=1}^{l}(1+(-tuv)^{j})^{b_{j}r}
\end{eqnarray*}
where $c_{[a,b]}$ is the size of the conjugacy class in $W_{Sp_{n}}$,
corresponding to $[a,b]\in\mathcal{B}_{n}$.
\end{theorem}
The same method allows to obtain explicit expressions for 
$\mu(\mathcal{X}_{\mathbb{Z}^{r}}^{0}G)$
in the case of other reductive $G$; the special orthogonal groups $SO_n$ will be
addressed in a future work.

\section{Generating functions for $E$-polynomials}

In this section we consider character varieties with arbitrarily bad
singularities. In this case, there are formidable difficulties in
computing the corresponding Poincar\'e polynomials in general, and previous
explicit methods have dealt with the $E$-polynomials for low dimensional
groups such as $SL_{2}$ and $SL_{3}$ (\cite{LMN,LM,BH}). 

By using the additive and multiplicative properties of $E$-polynomial,
for $G=GL_{n}$ we now address our new approach on $E$-polynomial
computations based on a stratification of $\mathcal{X}_{\Gamma}G$
that we term \emph{by partition type}, and which works for arbitrary
$\Gamma$.


\label{sec:stratification-partition} Using standard arguments in
GIT, any character variety admits a stratification by the dimension
of the stabilizer of a given representation. When $G$ is the general
linear group $GL_{n}$ (as well as the related groups $SL_{n}$ and
$PGL_{n}$), there is a more convenient refined stratification that
gives a lot of information on the corresponding character varieties
$\mathcal{X}_{\Gamma}G$ which we call \emph{stratification by partition
type}. 

\begin{definition}
Let $G=GL_{n}$ and $[k]\in\mathcal{P}_{n}$. We say that $\rho\in\mathcal{R}_{\Gamma}G=\hom(\Gamma,G)$
is $[k]$-polystable if $\rho$ is conjugated to $\bigoplus_{j=1}^{n}\rho_{j}$
where each $\rho_{j}$ is, in turn, a direct sum of $k_{j}>0$ \emph{irreducible}
representations of $\mathcal{R}_{\Gamma}(GL_{j})$, for $j=1,\cdots,n$
(by convention, if some $k_{j}=0$, then $\rho_{j}$ is not present
in the direct sum). 
\end{definition}

We denote $[k]$-polystable representations by $\mathcal{R}_{\Gamma}^{[k]}G$
and use similar terminology/notation for equivalence classes under
conjugation $\mathcal{X}_{\Gamma}^{[k]}G\subset\mathcal{X}_{\Gamma}G$.
It is to be noted that the trivial partition $[n]=[n^{1}]\in\mathcal{P}_{n}$ corresponds exactly to the irreducible (or stable) locus: 
$\mathcal{X}_{\Gamma}^{[n]}G=\mathcal{X}_{\Gamma}^{irr}G$. 

\begin{proposition}
\label{prop:locally-closed-GLn}Fix $n\in\mathbb{N}$, and let $G=GL_{n}$.
Then $\mathcal{X}_{\Gamma}G=\bigsqcup_{[k]\in\mathcal{P}_{n}}\mathcal{X}_{\Gamma}^{[k]}G,$
as a disjoint union of locally closed quasi-projective varieties.
\end{proposition}


The next result relates, by the plethystic exponential, the generating
functions of the $E$-polynomials $E(\mathcal{X}_{\Gamma}GL_{n})$
to the corresponding generating functions of the $E$-polynomials
of the irreducible character varieties $E(\mathcal{X}_{\Gamma}^{irr}GL_{n})$.

The \emph{plethystic exponential} of a formal power series $f(x,y,z)=\sum_{n\geq0}f_{n}(x,y)\,z^{n}\in\mathbb{Q}[x,y][[z]]$
is denoted by $\pexp(f)$, and defined formally (in terms of the usual exponential)
as $\pexp(f):=e^{\Psi(f)}\in\mathbb{Q}[x,y][[z]],$ where $\Psi$
acts on monomials as: $\Psi(x^{i}y^{j}z^{k})=\sum_{l\geq1}\frac{x^{li}y^{lj}z^{lk}}{l},$
where $(i,j,k)\in\mathbb{N}_{0}^{3}\setminus\{(0,0,0)\}$, and is $\mathbb{Q}$-linear
on $\mathbb{Q}[x,y][[z]]$. This exponential plays a prominent role
in the combinatorics of symmetric functions, and has applications
in counting of gauge invariant operators in supersymmetric quantum
field theories (see eg. \cite{FHH}).

\begin{theorem}
\label{thm:A-PExpB}Let $\Gamma$ be any finitely presented group.
Then: 
\[
\sum_{n\geq0}E(\mathcal{X}_{\Gamma}GL_{n};u,v)\,t^{n}=\pexp\left(\sum_{n\geq1}E(\mathcal{X}_{\Gamma}^{irr}GL_{n};u,v)\,t^{n}\right).
\]
\end{theorem}
The proofs of Theorem \ref{thm:A-PExpB} and Proposition \ref{prop:locally-closed-GLn}
are detailed in \cite{FNZ}; they allow to write explicit expressions for $E(\mathcal{X}_{\Gamma}GL_{n})$, for any group $\Gamma$, for which we have a formula for $E(\mathcal{X}_{\Gamma}^{irr}GL_{m})$, for all $m\leq n$, by a simple finite algorithm (and vice-versa). The formula
of Theorem \ref{thm:A-PExpB} generalizes a formula of \cite{MR} 
to an arbitrary group $\Gamma$,
even if the corresponding $GL_{n}$-character variety is not of polynomial type.

\section{The Free Group Case}

In this last section, we describe applications of the above methods
to the case of the free group of rank $r$, $\Gamma=F_{r}$; for simplicity
we adopt the notations $\mathcal{X}_{r}GL_{n}$, $\mathcal{X}_{r}SL_{n}$,
etc, for the corresponding character varieties. In \cite{MR}, it
was shown that $\mathcal{X}_{r}^{irr}GL_{n}$ and $\mathcal{X}_{r}GL_{n}$
are of polynomial type. Moreover, by counting points over finite fields
and using a theorem of Katz (\cite[Appendix]{HRV}), Mozgovoy and
Reineke found a formula for the $E$-polynomial of $\mathcal{X}_{r}^{irr}GL_{n}$
that can be written as follows (dropping the $x$ variable in $E(X;\,x)$,
and using $|[k]|:=k_{1}+\cdots+k_{d}$ for the length of a partition $[k]\in\mathcal{P}_{d}$).
\begin{proposition}
\cite{MR,FNZ}\label{prop:E-irr} For $r,n\geq2$, we have: 
\[
E(\mathcal{X}_{r}^{irr}GL_{n})=(x-1)\sum_{d|n}\frac{\mu(n/d)}{n/d}\,\sum_{[k]\in\mathcal{P}_{d}}\frac{(-1)^{|[k]|}}{|[k]|}\binom{|[k]|}{k_{1},\cdots,k_{d}}\prod_{j=1}^{d}b_{j}(x^{n/d})^{k_{j}}x^{\frac{n(r-1)k_{j}}{d}\binom{j}{2}},
\]
where $\mu$ is the M\"obius function, and the $b_{j}(x)$ are polynomials defined
by:\\[-3mm] 
\begin{equation}
(1+\sum_{n\geq1}b_{n}(x)\,t^{n})\left(1+\sum_{n\geq1}\big((x-1)(x^{2}-1)\ldots(x^{n}-1)\big)^{r-1}\,t^{n}\right)=1.\label{eq:F(t)}
\end{equation}

\end{proposition}
Using Propositions \ref{prop:locally-closed-GLn} and \ref{prop:E-irr}
and Theorem \ref{thm:A-PExpB}, we are able to write down very explicit
expressions for $E(\mathcal{X}_{r}^{[k]}GL_{n})$, the $E$-polynomials
of all polystable strata of $\mathcal{X}_{r}GL_{n}$ (see \cite[Secs. 5 and 6]{FNZ},
where we also compute $E(\mathcal{X}_{\Gamma}^{irr}GL_{n})$ for other
$\Gamma$ and low $n$).

We now provide a few lines on a forthcoming proof of the equality
between the $E$-polynomials of $\mathcal{X}_{r}SL_{n}$ and of $\mathcal{X}_{r}PGL_{n}$
for all $n\in\mathbb{N}$. This has been conjectured in Lawton-Mu\~noz
in \cite{LM}, who proved by explicit computation the cases $n=2$
and 3. 

In a analogous way as for $GL_{n}$ (see Section~\ref{sec:stratification-partition}),
we can define the $[k]$-polystable loci $\mathcal{X}_{r}^{[k]}SL_{n}$
and $\mathcal{X}_{r}^{[k]}PGL_{n}$ as follows. For a partition $[k]\in\mathcal{P}_{n}$,
the $[k]$-stratum of $\mathcal{X}_{r}SL_{n}$ is defined by restriction
of the corresponding one for $GL_{n}$: 
\[
\mathcal{X}_{r}^{[k]}SL_{n}:=\{\rho\in\mathcal{X}_{r}^{[k]}GL_{n}\,|\,\det\rho=1\},
\]
where the determinant of a representation is an element of $\mathcal{R}_{r}\mathbb{C}^{*}$.
By considering the action $\mathcal{R}_{r}\mathbb{C}^{*}\times\mathcal{X}_{r}GL_{n}\to\mathcal{X}_{r}GL_{n}$
given by multiplication of (conjugacy classes of) representations,
which is well defined on the GIT quotients and preserves the stratification
of $GL_{n}$, we can define 
\begin{equation}
\mathcal{X}_{r}^{[k]}PGL_{n}:=\mathcal{X}_{r}^{[k]}GL_{n}/\mathcal{R}_{r}\mathbb{C}^{*}=\mathcal{X}_{r}^{[k]}GL_{n}/(\mathbb{C}^{*})^{r}.\label{eq:GLn-PGLn}
\end{equation}

\begin{theorem}
\label{thm:SLn-PGLn}\cite{FNZ2} For the free group $F_{r}$, we
have the equalities: 
\begin{eqnarray*}
E(\mathcal{X}_{r}SL_{n}) & = & E(\mathcal{X}_{r}PGL_{n})=E(\mathcal{X}_{r}GL_{n})(x-1)^{-r}\\
E(\mathcal{X}_{r}^{[k]}SL_{n}) & = & E(\mathcal{X}_{r}^{[k]}PGL_{n})=E(\mathcal{X}_{r}^{[k]}GL_{n})(x-1)^{-r},
\end{eqnarray*}
for every $r,\;n$ and partition $[k]\in\mathcal{P}_{n}$.
\end{theorem}
The proof of Theorem \ref{thm:SLn-PGLn} uses geometric methods and
has two parts. The easy part is the relation between the $E$-polynomials
of $\mathcal{X}_{r}^{[k]}PGL_{n}$ and of $\mathcal{X}_{r}^{[k]}GL_{n}$,
which follows from the locally trivial (in the Zariski topology) fibration 
corresponding to the quotient \eqref{eq:GLn-PGLn}. The difficult part
is the relation between the strata $\mathcal{X}_{r}^{[k]}PGL_{n}$
and $\mathcal{X}_{r}^{[k]}SL_{n}$ which involves finite quotients: it
requires the proof of the triviality of the action of the center 
$\mathbb{Z}_{n}\subset SL_{n}$
on the cohomology (with compact support) of all the 
strata $\mathcal{X}_{r}^{[k]}SL_{n}$; for this
we use equivariant cohomology and a deformation retraction
between $\mathcal{X}_{r}^{irr}SL_{n}$ and the smooth part of the
semialgebraic set $\hom(F_r,SU(n))/SU(n)$ (see \cite{FL}).


\end{document}